\documentclass[12pt,reqno,english]{article}
\usepackage[T1]{fontenc}
\usepackage[latin9]{inputenc}
\usepackage{geometry}
\usepackage{color}
\usepackage{babel}
\usepackage{amsmath}
\usepackage{amsthm}
\usepackage{amssymb}
\usepackage{graphicx}
\PassOptionsToPackage{normalem}{ulem}
\usepackage{ulem}
\usepackage[unicode=true,pdfusetitle,
 bookmarks=true,bookmarksnumbered=false,bookmarksopen=false,
 breaklinks=false,pdfborder={0 0 1},backref=page,colorlinks=true]
 {hyperref}
\hypersetup{
 linkcolor=blue,citecolor=blue}

\makeatletter
\numberwithin{equation}{section}
\numberwithin{figure}{section}
\theoremstyle{plain}
\newtheorem{thm}{\protect\theoremname}[section]
\theoremstyle{plain}
\newtheorem{cor}[thm]{\protect\corollaryname}
\theoremstyle{remark}
\newtheorem{rem}[thm]{\protect\remarkname}
\theoremstyle{plain}
\newtheorem{prop}[thm]{\protect\propositionname}
\theoremstyle{plain}
\newtheorem{lem}[thm]{\protect\lemmaname}
\theoremstyle{plain}
\newtheorem{conjecture}[thm]{\protect\conjecturename}

\usepackage[leftcaption]{sidecap}

\sidecaptionvpos{figure}{c}

\newcommand{\FigBesBeg}[1][1.0]{%
 \let\MyFigure\figure
 \let\MyEndfigure\endfigure
 \renewenvironment{figure}[1]{\begin{SCfigure}[#1]##1}{\end{SCfigure}}}

\newcommand{\FigBesEnd}{%
 \let\figure\MyFigure
 \let\endfigure\MyEndfigure}

\usepackage{bbm}

\usepackage{bbold} % for the characteristic function symbol \mathbb{1}

\usepackage[all]{xy}

\theoremstyle{plain}
\newtheorem{claim}[thm]{\protect\claimname}

\makeatother

\providecommand{\conjecturename}{Conjecture}
\providecommand{\corollaryname}{Corollary}
\providecommand{\lemmaname}{Lemma}
\providecommand{\propositionname}{Proposition}
\providecommand{\remarkname}{Remark}
\providecommand{\theoremname}{Theorem}

\begin{document}
\global\long\def\F{\mathbb{\mathbb{\mathbf{F}}}}%
 
\global\long\def\defi{\stackrel{\mathrm{def}}{=}}%
 
\global\long\def\cogr{\mathrm{cogr}}%
 
\global\long\def\sub{S}%

\title{An extension of the cogrowth formula\\
to arbitrary subsets of the tree}
\author{Doron Puder}
\maketitle
\begin{abstract}
What is the probability that a random walk in the free group ends
in a proper power? Or in a primitive element? We present a formula
that computes the exponential decay rate of the probability that a
random walk on a regular tree ends in a given subset, in terms of
the exponential decay rate of the analogous probability of the non-backtracking
random walk. This generalizes the well-known cogrowth formula of Grigorchuk
\cite{Gri80}, Cohen \cite{cohen1982cogrowth} and Northshield \cite{Nor92}.
We also extend the formula to arbitrary subsets of the biregular tree.
\end{abstract}
\tableofcontents{}

\section{Introduction}

Let $d\ge3$ be an integer, let $\mathbb{T}_{d}$ denote the $d$-regular
tree with a root vertex $o$, and let $\sub\subseteq V\left(\mathbb{T}_{d}\right)$
be some subset of the vertices. Denote by $a_{t}\left(\sub\right)$
the number of non-backtracking walks of length $t$ in $\mathbb{T}_{d}$
that begin in $o$ and end inside $\sub$. Equivalently, $a_{t}\left(\sub\right)$
is the number of elements of $S$ at distance $t$ from $o$. Let
\begin{equation}
\alpha\left(\sub\right)\defi\limsup_{t\to\infty}a_{t}\left(\sub\right)^{1/t}\label{eq:def of alpha}
\end{equation}
denote the exponential growth rate of these numbers. Similarly, let
$b_{t}\left(\sub\right)$ denote the number of arbitrary, not necessarily
non-backtracking, walks of length $t$ that begin in $o$ and end
inside $\sub$, and denote 
\begin{equation}
\beta\left(\sub\right)\defi\limsup_{t\to\infty}b_{t}\left(\sub\right)^{1/t}.\label{eq:def of beta}
\end{equation}
Of course, $\frac{a_{t}\left(\sub\right)}{d\left(d-1\right)^{t-1}}$
is the probability that the length-$t$ non-backtracking simple random
walk on $\mathbb{T}_{d}$ starting at $o$ ends inside $\sub$, and
so $\frac{\alpha\left(\sub\right)}{d-1}$ is the exponential growth
rate (or, rather, decay rate) of this probability. Likewise, $\frac{\beta\left(\sub\right)}{d}$
is the exponential growth rate of the analogous probability for the
simple random walk.

For certain subsets $\sub$ it is easier to compute directly one of
these numbers, usually $\alpha\left(\sub\right)$, than it is to directly
compute the other. For example, if $d=2r$, then $\mathbb{T}_{d}$
is also the Cayley graph of $\F_{r}$, the free group of rank $r$,
with respect to some basis. Let $\sub$ denote the subset of squares
in $\F_{r}$. Every reduced word in $\F_{r}$ which is a square is
of the reduced form $uv^{2}u^{-1}$ where $u$ and $v$ are some (reduced)
words in $\F_{r}$. While there are no squares of odd length, for
every even $t$, the number of such words is, up to constants, $t\left(2r-1\right)^{t/2}$:
the factor $t$ is for the length $\ell\in\left[0,\frac{t}{2}\right]$
of $u$, which, in turn, determines the length of $v$ as $\frac{t}{2}-\ell$,
and then there are $\left(2r-1\right)^{\ell}$ options for $u$ and
$\left(2r-1\right)^{t/2-\ell}$ options for $v$ (again, up to constants).
So in this case, $\alpha\left(\sub\right)=\sqrt{2r-1}=\sqrt{d-1}$.
Does this number determine\footnote{A very similar argument shows that when $\sub$ is the subset of all
proper powers in $\F_{r}$ we also have $\alpha\left(\sub\right)=\sqrt{2r-1}=\sqrt{d-1}$.
We remark that the number $\beta\left(\sub\right)$ for this $\sub$
is used, for example, in the classical proof of Broder and Shamir
that random $d$-regular graphs are expanders \cite{BS87} (the argument
also appears in \cite[\S7.2]{HLW06}). In that work, $\beta\left(\sub\right)$
is bounded from above by $2\sqrt{d}$. As we note below in Corollary
\ref{cor:proper powers}, the precise value is $2\sqrt{d-1}$.} $\beta\left(\sub\right)$?

For certain sets $\sub$, it is known that $\beta\left(\sub\right)$
is a function of $\alpha\left(\sub\right)$. Let $p\colon\mathbb{T}_{d}\to G$
be a covering map onto a connected graph ($d$-regular, of course)
and $\sub=p^{-1}\left(v\right)$ be the fiber above some vertex $v$
of $G$. In this case, $\alpha\left(\sub\right)$ is the \emph{cogrowth}
of $G$, also denoted $\mathrm{cogr}\left(G\right)$. This is the
exponential growth rate of the number of non-backtracking walks from
$v$ to itself (or, equivalently, the number of non-backtracking walks
from $u$ to $v$ for any two vertices $u,v\in V\left(G\right)$:
the exponential growth rate does not depend on $u$ and $v$). This
term was first defined for Cayley graphs by Grigorchuk \cite{Gri80}
and, independently, by Cohen \cite{cohen1982cogrowth}, and was later
extended to arbitrary regular graphs by Northshield \cite{Nor92}.
In this case, the number $\beta\left(\sub\right)$ is the exponential
growth rate of the number of not-necessarily-non-backtracking walks
from $v$ to itself (or, again, from $u$ to $v$ for arbitrary $u,v,\in V\left(G\right)$),
and is known to be equal to the spectral radius of the adjacency operator
of $G$, which we denote here $\rho\left(G\right)$. For this type
of sets, a formula for $\beta\left(\sub\right)=\rho\left(G\right)$
in terms of $\alpha\left(\sub\right)=\mathrm{cogr}\left(G\right)$
is known as \emph{the cogrowth formula}:
\begin{thm}[Classical cogrowth formula]
\label{thm:cogrowth-formula}\cite{Gri80,cohen1982cogrowth,Nor92}
Let $G$ be a connected $d$-regular graph ($d\ge3$). Then 
\[
\rho\left(G\right)=\begin{cases}
2\sqrt{d-1} & \mathrm{if}~\mathrm{cogr}\left(G\right)\in\left[0,\sqrt{d-1}\right]\\
\mathrm{cogr}\left(G\right)+\frac{d-1}{\mathrm{cogr}\left(G\right)} & \mathrm{if}~\mathrm{cogr}\left(G\right)\geq\left[\sqrt{d-1},d-1\right]
\end{cases}.
\]
\end{thm}

Note that for every $G$, $\mathrm{cogr}\left(G\right)\in\left[0,d-1\right]$
and $\rho\left(G\right)\in\left[2\sqrt{d-1},d\right]$. In addition,
$\mathrm{cogr}\left(G\right)=0$ if and only if $G=\mathbb{T}_{d}$
is the tree, and otherwise $\mathrm{cogr}\left(G\right)\ge1$. For
non-tree Cayley graphs, $\cogr\left(G\right)\ge\sqrt{d-1}$ \cite[Thm.~1]{cohen1982cogrowth}.
Not surprisingly, a very similar formula connects the eigenvalues
of the adjacency matrix of a finite $d$-regular graph and the eigenvalues
of the non-backtracking Hashimoto matrix: this formula for the eigenvalues
follows from the Ihara-Bass formula (see, e.g., \cite{kotani2000zeta}
or \cite[\S2]{friedman2023note}), although we are not aware of a
direct derivation of one from the other. \\

In this paper we show that the same formula holds for arbitrary subsets
of the tree. Denote by $g_{d}\colon\mathbb{R}_{\ge0}\to[2\sqrt{d-1},\infty)$
the map appearing in the cogrowth formula, namely,
\begin{equation}
g_{d}\left(\alpha\right)=\begin{cases}
2\sqrt{d-1} & \mathrm{if}~\alpha\le\sqrt{d-1},\\
\alpha+\frac{d-1}{\alpha} & \mathrm{if}~\alpha\ge\sqrt{d-1}.
\end{cases}\label{eq:def of g}
\end{equation}

\begin{thm}[Extended cogrowth formula]
\label{thm:cogrowth for subsets} Let $\emptyset\ne\sub\subseteq V\left(\mathbb{T}_{d}\right)$
be an arbitrary non-empty subset of the $d$-regular tree. Then,
\begin{enumerate}
\item Let $\alpha\left(\sub\right)=\limsup_{t\to\infty}a_{t}\left(\sub\right)^{1/t}$
and $\beta\left(\sub\right)=\limsup_{t\to\infty}b_{t}\left(\sub\right)^{1/t}$
be defined as above. Then 
\[
\beta\left(\sub\right)=g_{d}\left(\alpha\left(S\right)\right).
\]
\item The same formula holds for ordinary limits: denote by $S^{\mathrm{even}}$
($S^{\mathrm{odd}}$) the subset of vertices of $S$ at even (respectively,
odd) distance from $o$.
\begin{enumerate}
\item If $\tilde{\alpha}\left(S\right)\defi\lim_{t\to\infty}a_{t}\left(\sub\right)^{1/t}$
exists and $S^{\mathrm{even}},S^{\mathrm{odd}}\ne\emptyset$, then
$\tilde{\beta}\left(\sub\right)\defi\lim_{t\to\infty}b_{t}\left(\sub\right)^{1/t}$
exists, and $\tilde{\beta}\left(S\right)=g_{d}\left(\tilde{\alpha}\left(S\right)\right)$.
\item If 
\[
\tilde{\alpha}^{\mathrm{even}}\left(S\right)\defi\lim_{t\to\infty}a_{2t}\left(S\right)^{1/(2t)}
\]
exists and $S^{\mathrm{even}}\ne\emptyset$, then $\tilde{\beta}^{\mathrm{even}}\left(S\right)\defi\lim_{t\to\infty}b_{2t}\left(S\right)^{1/(2t)}$
exists and $\tilde{\beta}^{\mathrm{even}}\left(S\right)=g_{d}\left(\tilde{\alpha}^{\mathrm{even}}\left(S\right)\right)$.
\item If 
\[
\tilde{\alpha}^{\mathrm{odd}}\left(S\right)\defi\lim_{t\to\infty}a_{2t+1}\left(S\right)^{1/(2t+1)}
\]
exists and $S^{\mathrm{odd}}\ne\emptyset$, then $\tilde{\beta}^{\mathrm{odd}}\left(S\right)\defi\lim_{t\to\infty}b_{2t+1}\left(S\right)^{1/(2t+1)}$
exists and $\tilde{\beta}^{\mathrm{odd}}\left(S\right)=g_{d}\left(\tilde{\alpha}^{\mathrm{odd}}\left(S\right)\right)$.
\end{enumerate}
\end{enumerate}
\end{thm}

Theorem \ref{thm:cogrowth for functions} below extends the cogrowth
formula further to arbitrary non-negative functions defined on the
vertices of the tree. Most interestingly, this allows for exponential
growth rates that are arbitrarily large, unlike the original formula
which is relevant only to $\alpha\left(S\right)\le d-1$ and $\beta\left(S\right)\le d$.
We show that the same formula ($g_{d}$ as defined in \eqref{eq:def of g})
still applies in this larger domain.

Since it was first discovered, quite a few proofs have been found
for the original cogrowth formula (Theorem \ref{thm:cogrowth-formula}):
these include \cite{Gri80,cohen1982cogrowth,szwarc1989short,Nor92,Bar99,LyonsPeres2017},
the simplest of which is probably the one in \cite[Thm. 6.10]{LyonsPeres2017}.
Our proof of the more general Theorems \ref{thm:cogrowth for subsets}
and \ref{thm:cogrowth for functions} gives yet another proof also
for the original formula. \\

Let us mention some special cases of Theorem \ref{thm:cogrowth for subsets}.
Above we explained why $\alpha\left(\sub\right)=\sqrt{d-1}$ for the
set $\sub$ of squares in $\F_{r}$ ($r\ge2$), where $d=2r$. It
follows immediately that $\beta\left(S\right)=2\sqrt{d-1}$. Moreover,
\begin{cor}
\label{cor:proper powers}Let $d=2r$, let $\mathbb{T}_{d}$ be the
Cayley graph of the free group $\F_{r}$ ($r\ge2$) with respect to
some basis, and let $\sub$ be the subset of proper powers. Then $\beta\left(\sub\right)=2\sqrt{d-1}$.
\end{cor}

\begin{proof}
This follows from Theorem \ref{thm:cogrowth for subsets} using the
fact that $\alpha\left(\sub\right)=\sqrt{d-1}$. The latter fact is
an easy adaptation of the argument for squares above. It also appears,
for example, as the special case $m=1$ in \cite[Thm.~8.2]{Pud15}.
\end{proof}
\begin{cor}
\label{cor:primitives}Let $d=2r$, let $\mathbb{T}_{d}$ be the Cayley
graph of the free group $\F_{r}$ ($r\ge2$) with respect to some
basis, and let ${\cal P}_{r}$ be the subset of primitive words (words
belonging to some free basis of $\F_{r}$). Then 
\[
\beta\left({\cal P}_{r}\right)=\begin{cases}
2\sqrt{3} & \mathrm{if}~r=2,\\
2r-3+\frac{2r-1}{2r-3} & \mathrm{if}~r\ge3.
\end{cases}
\]
\end{cor}

\begin{proof}
This follows from the values of $\alpha\left({\cal P}_{r}\right)$
and Theorem \ref{thm:cogrowth for subsets}. That $\alpha\left({\cal P}_{r}\right)=\sqrt{3}$
in $\F_{2}$ is a result of Rivin \cite{rivin2004remark}. That $\alpha\left({\cal P}_{r}\right)=2r-3$
in $\F_{r}$, $r\ge3$, is the main result of \cite{PW14}.
\end{proof}
The original cogrowth formula (Theorem \ref{thm:cogrowth-formula})
can be also seen as the special case of the following corollary in
which the subset is $S=\left\{ o\right\} \subseteq V\left(G\right)$.
\begin{cor}
\label{cor:arbitrary regular graphs}The formulas of Theorem \ref{thm:cogrowth for subsets}
hold for subsets of the vertices of an arbitrary connected $d$-regular
graph $G$. 

For example, if $S\subseteq V\left(G\right)$, $o\in V\left(G\right)$
an arbitrary vertex, $a_{t}\left(S\right)$ is the number of non-backtracking
length-$t$ walks from $o$ into $S$, $b_{t}\left(S\right)$ is the
number of arbitrary length-$t$ walks from $o$ into $S$ and $\alpha\left(S\right)$
and $\beta\left(S\right)$ are their respective exponential growth
rate, then $\beta\left(S\right)=g_{d}\left(\alpha\left(S\right)\right)$
with $g_{d}$ as defined in \eqref{eq:def of g}.
\end{cor}

\begin{proof}
Pull back $S$ to the covering $d$-regular tree of $G$.
\end{proof}
\begin{rem}
\label{rem:application for primitivity rank}Another application of
the extended cogrowth formula lies in our original motivation for
this formula in \cite{Pud15}: the formula was used to compute (in
{[}ibid, Thm.~8.2{]}) the precise exponential growth rate of the
set of not-necessarily-reduced words in a free group with a given
'primitivity rank' \cite[Def.~1.7]{Pud14a}. These precise estimates
were recently used in the exciting new work by Chen, Garza-Vargas,
Tropp and van Handel \cite[Thm.~3.6]{chen2024new}.
\end{rem}

\subsection*{The biregular case}

There is a further extension of the cogrowth formula to arbitrary
subsets of the $\left(c,d\right)$-biregular tree. This is the tree
with vertices of degrees $c$ and $d$, where every vertex of degree
$d$ has all its neighbors of degree $c$, and vice versa. 
\begin{thm}[Biregular cogrowth formula]
\label{thm:cogrowth biregular} Let $\mathbb{T}_{c,d}$ be the $\left(c,d\right)$-biregular
tree rooted at a vertex $o$, with $2\le c<d$, and $\emptyset\ne\sub\subseteq V\left(\mathbb{T}_{c,d}\right)$
an arbitrary non-empty subset of its vertices. Let $\alpha\left(\sub\right)$
and $\beta\left(\sub\right)$ be defined as in \eqref{eq:def of alpha}
and \eqref{eq:def of beta}. Then 
\[
\beta\left(\sub\right)=\begin{cases}
\sqrt{c-1}+\sqrt{d-1} & \mathrm{if}~\alpha\left(\sub\right)\le\left[\left(c-1\right)\left(d-1\right)\right]^{1/4},\\
\sqrt{\alpha\left(S\right)^{2}+c+d-2+\frac{\left(c-1\right)\left(d-1\right)}{\alpha\left(S\right)^{2}}} & \mathrm{if}~\alpha\left(\sub\right)\ge\left[\left(c-1\right)\left(d-1\right)\right]^{1/4}.
\end{cases}
\]
\end{thm}

\subsection*{Outline of the paper}

Section \ref{sec:regular} proves our main result -- Theorem \ref{thm:cogrowth for subsets}
and the more general version applying to function and stated as Theorem
\ref{thm:cogrowth for functions}. In Section \ref{sec:biregular}
we prove the biregular version -- Theorem \ref{thm:cogrowth biregular}
-- using a simpler argument, which does not apply to arbitrary non-negative
functions, but contains the case of subsets in the regular case.

\subsection*{Acknowledgments}

This paper was originally planned to appear much earlier. The main
result of this paper was first announced in \cite[Thm.~4.4]{Pud15}
with the application mentioned in Remark \ref{rem:application for primitivity rank},
which was then used to analyze the adjacency operator in random regular
graphs. But the writing was not completed back then, partly because
the original motivation -- its usage in the proof of the main result
of \cite{Pud15} -- became obsolete: a short time afterwards we realized
that analyzing the non-backtracking spectrum directly is more fruitful
-- see \cite{friedman2023note} and especially Remark 3.4 therein.
We thank Ramon van Handel for encouraging us to finally complete the
writing of the current paper (his interest stemmed from his new work
with co-authors \cite{chen2024new} where they used the application
from Remark \ref{rem:application for primitivity rank}). We thank
Asaf Nachmias and Ori Parzanchevski for beneficial conversations.
This work was supported by the European Research Council (ERC) under
the European Union\textquoteright s Horizon 2020 research and innovation
programme (grant agreement No 850956), by the Israel Science Foundation,
ISF grants 1140/23, as well as by the National Science Foundation
under Grant No. DMS-1926686. 

\section{The cogrowth formula for arbitrary subsets of the regular tree\label{sec:regular}}

Fix $d\in\mathbb{Z}$, $d\ge3$. In this section we prove our main
result, Theorem \ref{thm:cogrowth for subsets}, that for arbitrary
$\emptyset\ne\sub\subseteq V\left(\mathbb{T}_{d}\right)$ we have
$\beta\left(\sub\right)=g_{d}\left(\alpha\left(\sub\right)\right)$.
In fact, we prove a slightly more general result. Let $0\ne f\colon V\left(\mathbb{T}_{d}\right)\to\mathbb{R}_{\ge0}$
be a non-negative function defined on the vertices of the $d$-regular
tree which is not identically zero. We define 
\begin{align*}
a_{t}\left(f\right)=\sum_{\substack{p:\,\,\mathrm{a\,non\text{-}backtracking\,walk}\\
\mathrm{from}\,o\,\mathrm{of\,length}\,t
}
}f\left(\mathrm{end}\left(p\right)\right)~~~~~~~~~~ &  & \alpha\left(f\right)=\limsup_{t\to\infty}a_{t}\left(f\right)^{1/t}\\
b_{t}\left(f\right)=\sum_{\substack{p:\,\,\mathrm{a\,walk\,from\,}o\\
\mathrm{of\,length}\,t
}
}f\left(\mathrm{end}\left(p\right)\right)~~\,\,\,\,\,\,\,\,\,\,\,\,\,\,\,\,\,\,~\,\,\,~~~~~~~ &  & \,\beta\left(f\right)=\limsup_{t\to\infty}b_{t}\left(f\right)^{1/t}.
\end{align*}
In principle, the values of $\alpha\left(f\right)$ and $\beta\left(f\right)$
may be $\infty$ if \textbf{$a_{t}\left(f\right)$ }(respectively
$b_{t}\left(f\right)$) grows super-exponentially, but as $b_{t}\left(f\right)\ge a_{t}\left(f\right)$
by definition, we have $\alpha\left(f\right)=\infty\Longrightarrow\beta\left(f\right)=\infty$.
We may thus assume that $\alpha\left(f\right)<\infty$. Of course,
if $\mathbb{1}_{\sub}$ is the characteristic function of the subset
$\sub$, then $a_{t}\left(\mathbb{1}_{\sub}\right)=a_{t}\left(\sub\right)$,
$\alpha\left(\mathbb{1}_{\sub}\right)=\alpha\left(\sub\right)$ and
so on. So the following is a generalization of Theorem \ref{thm:cogrowth for subsets}.
\begin{thm}
\label{thm:cogrowth for functions} Fix $d\ge3$ and let $0\ne f\colon V\left(\mathbb{T}_{d}\right)\to\mathbb{R}_{\ge0}$
be a non-negative function on the vertices of the $d$-regular tree,
not identically zero, with $\alpha\left(f\right)<\infty$. Then,
\begin{enumerate}
\item \label{enu:limsup}$\beta\left(f\right)=g_{d}\left(\alpha\left(f\right)\right)$.
\item \label{enu:ordinary limit}The same formula holds for ordinary limits:
denote by $\mathbb{T}_{d}^{\mathrm{even}}$ ($\mathbb{T}_{d}^{\mathrm{odd}}$)
the subset of vertices of $\mathbb{T}_{d}$ at even (respectively,
odd) distance from $o$.
\begin{enumerate}
\item \label{enu:limit-odd-and-even}If $\tilde{\alpha}\left(f\right)\defi\lim_{t\to\infty}a_{t}\left(f\right)^{1/t}$
exists and $f|_{\mathbb{T}_{d}^{\mathrm{even}}},f|_{\mathbb{T}_{d}^{\mathrm{odd}}}\ne0$,
then $\tilde{\beta}\left(f\right)\defi\lim_{t\to\infty}b_{t}\left(f\right)^{1/t}$
exists, and $\tilde{\beta}\left(f\right)=g_{d}\left(\tilde{\alpha}\left(f\right)\right)$.
\item \label{enu:limit-even}If $\tilde{\alpha}^{\mathrm{even}}\left(f\right)\defi\lim_{t\to\infty}a_{2t}\left(f\right)^{1/(2t)}$
exists and $f|_{\mathbb{T}_{d}^{\mathrm{even}}}\ne0$, then $\tilde{\beta}^{\mathrm{even}}\left(f\right)\defi\lim_{t\to\infty}b_{2t}\left(f\right)^{1/(2t)}$
exists and $\tilde{\beta}^{\mathrm{even}}\left(f\right)=g_{d}\left(\tilde{\alpha}^{\mathrm{even}}\left(f\right)\right)$.
\item \label{enu:limit-odd}If $\tilde{\alpha}^{\mathrm{odd}}\left(f\right)\defi\lim_{t\to\infty}a_{2t+1}\left(f\right)^{1/(2t+1)}$
exists and $f|_{\mathbb{T}_{d}^{\mathrm{odd}}}\ne0$, then\linebreak{}
$\tilde{\beta}^{\mathrm{odd}}\left(f\right)\defi\lim_{t\to\infty}b_{2t+1}\left(f\right)^{1/(2t+1)}$
exists and $\tilde{\beta}^{\mathrm{odd}}\left(f\right)=g_{d}\left(\tilde{\alpha}^{\mathrm{odd}}\left(f\right)\right)$.
\end{enumerate}
\end{enumerate}
\end{thm}

The more general version of the theorem is useful in some natural
scenarios. For example, in \cite[Cor.~4.5]{Pud15}, we actually count
not just words of a given primitivity rank $m$, but the sum, over
all such words, of the number of ``critical subgroups'' of the word
(see \cite[Def.~2.2]{Pud15}). To illustrate, the special case of
$m=1$ is that of words that are proper powers, and then the corresponding
function $f\colon\F_{r}\to\mathbb{R}_{\ge0}$ assigns to every word
$w=u^{q}$, where $q\in\mathbb{Z}_{\ge1}$ and $u\in\F_{r}$ is a
non-power, the number $\frak{d}\left(q\right)-1$, where $\frak{d}\left(q\right)$
is the number of positive divisors of $q$.

\subsubsection*{The main idea of the proof}

\begin{figure}
\includegraphics[scale=0.3]{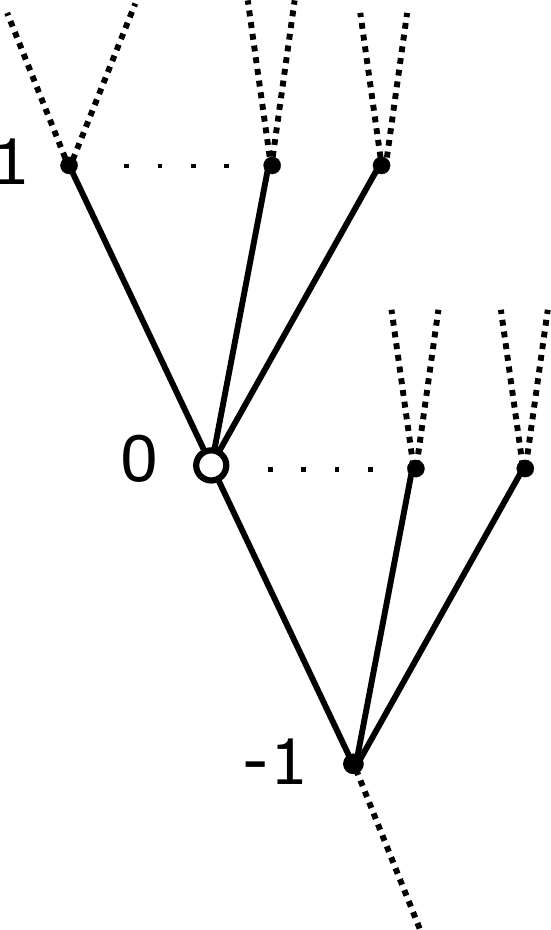}
\centering{}\caption{\label{fig:graded tree}$\mathbb{T}_{d}$ as a graded tree, with the
root $o$ marked with a non-filled circle; the vertices of the tree
are divided into levels in $\mathbb{Z}$}
\end{figure}

The main idea of the proof is to think of the $d$-regular tree $\mathbb{T}_{d}$
as a graded tree. Here, every vertex has a certain grade, or level,
$\text{\ensuremath{\ell\in}}\mathbb{Z}$, with one parent at level
$\ell-1$ and $d-1$ children at level $\ell+1$. The basepoint $o$
is at level $0$ -- see Figure \ref{fig:graded tree}. The number
of walks of length $t$ in the graded tree terminating at level $r$
is precisely $\binom{t}{y}\left(d-1\right)^{y}$, where $y=\frac{t+r}{2}$
is the number of steps going up the tree. We show that the latter
number approximates quite well the number of walks of length $t$
ending at distance $r$ from the root. More accurately, this is the
correct number up to a sub-exponential factor, so we may use it for
the sake of exponential growth rate.

For example, if $a_{t}\left(f\right)=\alpha^{t}$ for every $t$,
we thus obtain, up to sub-exponential factors, that
\begin{eqnarray}
b_{t}\left(f\right) & \approx & \sum_{r\in\mathbb{Z}\cap\left[0,t\right]\colon r\equiv t\pmod{2}}\binom{t}{\frac{r+t}{2}}\left(d-1\right)^{\left(r+t\right)/2}\frac{a_{r}\left(f\right)}{d\left(d-1\right)^{r-1}}\nonumber \\
 & \approx & \sum_{y=\left\lceil \frac{t}{2}\right\rceil }^{t}\binom{t}{y}\left(d-1\right)^{y}\left(\frac{\alpha}{d-1}\right)^{r}\nonumber \\
 & = & \sum_{y=\left\lceil \frac{t}{2}\right\rceil }^{t}\binom{t}{y}\left(\frac{d-1}{\alpha}\right)^{t-y}\alpha^{y}.\label{eq:half the terms}
\end{eqnarray}
The summation includes half of the terms of $\sum_{y=0}^{t}\binom{t}{y}\left(\frac{d-1}{\alpha}\right)^{t-y}\alpha^{y}=\left(\frac{d-1}{\alpha}+\alpha\right)^{t}$.
If we assume that $\alpha\ge\frac{d-1}{\alpha}$, or $\alpha\ge\sqrt{d-1}$,
the terms in \eqref{eq:half the terms} constitute more than half
of the weight, so up to a sub-exponential factor, $b_{t}\left(f\right)$
is $g_{d}\left(\alpha\right)^{t}$. For the case $\alpha\le\sqrt{d-1}$,
we use the argument in \eqref{eq:half the terms} to get an upper
bound $\beta\left(f\right)\le2\sqrt{d-1}$, and the fact that the
spectral radius of the adjacency operator in $\mathbb{T}_{d}$ is
$2\sqrt{d-1}$ as a lower bound. Below we make all these arguments
precise. 

\subsubsection*{The proof in detail}
\begin{prop}
\label{prop:approximate number of walks ending at distance r}Let
$t,y,r\in\mathbb{Z}_{\ge0}$ with $\frac{t}{2}\le y\le t$ and $r=2y-t$.
Denote by $w_{t}\left(o,r\right)$ the number of walks of length $t$
in $\mathbb{T}_{d}$ starting at $o$ and ending at distance $r$
from $o$. There exists a sub-exponential function $h\colon\mathbb{Z}_{\ge0}\to\mathbb{R}_{\ge0}$
(namely, $\limsup_{t\to\infty}h\left(t\right)^{1/t}=1$) so that 
\begin{equation}
\frac{1}{h\left(t\right)}\le\frac{w_{t}\left(o,r\right)}{\binom{t}{y}\left(d-1\right)^{y}}\le h\left(t\right).\label{eq:sub-exponential approx}
\end{equation}
\end{prop}

Note that $h\left(t\right)$ is independent of $y$ (and $r$): it
applies to every value of $y\in\left[\frac{t}{2},t\right]$. 

The proof of Proposition \ref{prop:approximate number of walks ending at distance r}
is given by the following Lemmas \ref{lem:positive and non-negative walks},
\ref{lem:lower bound}, \ref{lem:upper-bound, constant r} and \ref{lem:upper-bound, non-constant r}.
Our analysis uses the structure of a graded tree we impose on $\mathbb{T}_{d}$,
as in Figure \ref{fig:graded tree}. The term $\binom{t}{y}\left(d-1\right)^{y}$
is then precisely the number of walks of length $t$ starting at $o$
and ending at level $r=2y-t$ or, equivalently, with precisely $y$
steps upward. We call a walk in the graded tree \textbf{positive}
if except for the initial position it only visits positive levels.
We call it \textbf{non-negative} if it never visits negative levels.
Note that in a non-negative walk, the level of the ending vertex is
equal to its distance from the root $o$.
\begin{lem}
\label{lem:positive and non-negative walks}For every $t,y\in\mathbb{Z}_{\ge0}$
with $\frac{t}{2}\le y\le t$, the number of positive walks of length
$t$ and $y$ steps upward is precisely $\frac{2y-t}{t}\binom{t}{y}\left(d-1\right)^{y}$,
and the number of non-negative walks is precisely $\frac{2y-t+1}{y+1}\binom{t}{y}\left(d-1\right)^{y}$.
\end{lem}

\begin{proof}
For every length-$t$ walk corresponds a sequence of $t$ arrows $\nearrow$
and $\searrow$, describing whether any given step goes upward or
downward, respectively. The total number of walks with a prescribed
sequence with $y$ upward arrows and $t-y$ downward arrows is $\left(d-1\right)^{y}$
-- it does not depend on the locations of the $\nearrow$'s in the
sequence. So it is enough to show that $\frac{2y-t}{t}\binom{t}{y}$
of the $\binom{t}{y}$ possible sequences are positive. The argument
imitates one of the classical proofs for the fact that Catalan numbers
count Dyck paths. All such sequences end at level $2y-t$. For every
sequence with first step upward which is not positive, if one ``reflects''
it after the first positive time it reaches level $0$ and swaps every
$\nearrow$ with a $\searrow$ and vice versa, the sequence terminates
at height $t-2y$. This gives a bijection between sequences of length
$t$ which are not positive and with exactly $y$ steps upward including
the first step, and sequences of length $t$ with exactly $t-y$ steps
upward including the first step. The size of the latter set is $\binom{t-1}{t-y-1}$.
The total number of sequences of length $t$ with $y$ steps upward
including the first step is $\binom{t-1}{y-1}$. So the number of
positive sequences is 
\[
\binom{t-1}{y-1}-\binom{t-1}{t-y-1}=\binom{t-1}{y-1}-\binom{t-1}{y}=\frac{2y-t}{t}\binom{t}{y}.
\]
Similarly, every sequence that is not non-negative reaches level $-1$
at some point. If we reflect it after its first visit to level $-1$,
it ends at level $t-2y-2$. This gives a bijection between sequences
with $y$ steps forward which are not non-negative, and sequences
with exactly $t-y-1$ steps forward. The total number of non-negative
sequences is, therefore,
\[
\binom{t}{y}-\binom{t}{t-y-1}=\binom{t}{y}-\binom{t}{y+1}=\frac{2y-t+1}{y+1}\binom{t}{y}.
\]
\end{proof}
\begin{lem}
\label{lem:lower bound}In the notation of Proposition \ref{prop:approximate number of walks ending at distance r},
for every $t,y,r\in\mathbb{Z}_{\ge0}$ with $\frac{t}{2}\le y\le t$
and $r=2y-t$ we have
\[
\frac{1}{t}\le\frac{w_{t}\left(o,r\right)}{\binom{t}{y}\left(d-1\right)^{y}}.
\]
\end{lem}

\begin{proof}
Every non-negative walk of length $t$ and precisely $y$ steps upwards
ends at distance $2y-t$ from $o$. So if $r=2y-t$ we get that $w_{t}\left(o,r\right)$
is at least the number of such non-negative walks, and by Lemma \ref{lem:positive and non-negative walks},
\[
\frac{w_{t}\left(o,r\right)}{\binom{t}{y}\left(d-1\right)^{y}}\ge\frac{2y-t+1}{y+1}\ge\frac{1}{t}.
\]
\end{proof}
\begin{lem}
\label{lem:upper-bound, constant r}The upper bound in Proposition
\ref{prop:approximate number of walks ending at distance r} holds
for every constant $r_{0}\in\mathbb{Z}_{\ge0}$. Namely,
\[
\limsup_{t\to\infty}\left(\frac{w_{t}\left(o,r_{0}\right)}{\binom{t}{y}\left(d-1\right)^{y}}\right)^{1/t}\le1,
\]
where the limsup is over $t$'s with the same parity as $r_{0}$,
and $y=\frac{t+r_{0}}{2}$.
\end{lem}

\begin{proof}
This follows immediately from the fact that the spectral radius of
$\mathbb{T}_{d}$ is $2\sqrt{d-1}$. Indeed, by abuse of notation
let $w_{t}\left(u,v\right)$ denote the number of walks of length
$t$ from $u$ to $v$, where $u,v\in V\left(\mathbb{T}_{d}\right)$.
It is well-known that in $\mathbb{T}_{d}$, as in every locally finite
connected graph, the spectral radius of the adjacency operator is
equal to the ordinary limit $\lim_{t\to\infty}w_{2t}\left(u,v\right)^{1/\left(2t\right)}$
for arbitrary vertices $u$ and $v$ at even distance from each other,
or to $\lim_{t\to\infty}w_{2t+1}\left(u,v\right)^{1/\left(2t+1\right)}$
for arbitrary vertices at odd distance from one another. See, for
instance, \cite[page 191]{GZ99}. There are $d\left(d-1\right)^{r_{0}-1}$
vertices at distance $r_{0}$ from $o$ (or just one if $r_{0}=0$),
but a constant does not affect the exponential growth rate, and we
conclude that $\lim_{t\to\infty}w_{t}\left(o,r\right)^{1/t}=2\sqrt{d-1}$,
where the limit is over $t$'s with the same parity as $r_{0}$. But
a straightforward application of the Stirling approximation formula
yields that $\lim_{t\to\infty}\left[\binom{t}{y}\left(d-1\right)^{y}\right]^{1/t}=2\sqrt{d-1}$
where, again, the values of $t$ here have the same parity as $r_{0}$
and $y=\frac{t+r_{0}}{2}$.
\end{proof}
\begin{rem}
It seems, in fact, that for every constant $r_{0}$ there is a constant
$c>0$ so that in the notation of Lemma \ref{lem:upper-bound, constant r}
we have $\frac{w_{t}\left(o,r_{0}\right)}{\binom{t}{y}\left(d-1\right)^{y}}\le c$.
However, we do not need this stronger result here. 
\end{rem}

\begin{lem}
\label{lem:upper-bound, non-constant r}There is a constant $c>0$
so that in the notation of Proposition \ref{prop:approximate number of walks ending at distance r},
for every $t,y,r\in\mathbb{Z}_{\ge0}$ with $\frac{t}{2}\le y\le t$
and $r=2y-t\ge3$, we have
\[
\frac{w_{t}\left(o,r\right)}{\binom{t}{y}\left(d-1\right)^{y}}\le c.
\]
\end{lem}

\begin{proof}
Here we do \emph{not} use the grading of $\mathbb{T}_{d}$. The walks
of length $t$ in $\mathbb{T}_{d}$ reaching distance $r=2y-t$ can
be parameterized by sequences of $\nearrow$'s and $\searrow$'s,
expressing steps away from the root or steps towards the root, respectively.
Every sequence has precisely $y$ $\nearrow$'s. These sequences remain
``non-negative'': the distance to the root is always non-negative.
By Lemma \ref{lem:positive and non-negative walks}, the number of
such sequences is $\frac{2y-t+1}{y+1}\binom{t}{y}$. To count the
number of walks corresponding to a given sequence, note that there
are $d$ possible steps for every $\nearrow$ when the sequence is
at the root, $d-1$ possible steps for every other $\nearrow$, and
one option for every $\searrow$. Denote by $p_{i}=p_{i}(d,t,y)$,
$i=1,2,3,\ldots$, the proportion of sequences with precisely $i$
visits to the root (every sequence begins at the root, so the number
of visits is at least one). We get
\begin{equation}
w_{t}\left(o,r\right)=\frac{2y-t+1}{y+1}\binom{t}{y}\left(d-1\right)^{y}\left[p_{1}\left(\frac{d}{d-1}\right)+p_{2}\left(\frac{d}{d-1}\right)^{2}+p_{3}\left(\frac{d}{d-1}\right)^{3}+\ldots\right].\label{eq:precise formula for w_t(o,r)-1}
\end{equation}
We proceed by bounding from above the $p_{i}$'s. For every $i=1,2,3,\ldots$
let $q_{i}$ denote the proportion of sequences (of length $t$, with
$y$ $\nearrow$'s, and non-negative) visiting the root \emph{at least}
$i$ times. In particular, $q_{1}=1$. We bound $p_{i}$ from above
by $q_{i}$. Notice that $q_{2}$ is the proportion of sequences which
are not positive among all non-negative sequences, which, by Lemma
\ref{lem:positive and non-negative walks}, is
\begin{equation}
q_{2}=1-\frac{\frac{2y-t}{t}}{\frac{2y-t+1}{y+1}}=1-\frac{y+1}{t}\cdot\frac{2y-t}{2y-t+1}=1-\frac{y+1}{t}\cdot\frac{r}{r+1}.\label{eq:q_2}
\end{equation}
By assumption, $r\ge3$ and $\frac{y+1}{t}>\frac{y}{t}>\frac{1}{2}$,
so
\begin{equation}
q_{2}<1-\frac{1}{2}\cdot\frac{3}{4}=\frac{5}{8}.\label{eq:q_2 upper bound}
\end{equation}
We claim that $\frac{q_{3}}{q_{2}}\le q_{2}$. This should be intuitively
obvious: at the second visit of the root, the sequence needs to ``advance
faster'', on average, to reach distance $r$ at the end, than it
had at time zero. So we expect the proportion of sequences visiting
the root a third time among those visiting at least twice, to be smaller
than $q_{2}$. The formal argument is as follows. For $\tau=2,4,\ldots$,
denote\textbf{ }by $q_{2,\tau}$ the proportion of sequences with
second visit to the root at time $\tau$, and likewise by $q_{3,\tau}$
the proportion of sequences visiting the root at least thrice and
with \emph{second} visit to the root at time $\tau$. So 
\[
\frac{q_{3}}{q_{2}}=\frac{q_{3,2}+q_{3,4}+\ldots+q_{3,t-r}}{q_{2,2}+q_{2,4}+\ldots+q_{2,t-r}},
\]
which is a weighted average of $\frac{q_{3,\tau}}{q_{2,\tau}}$ for
$\tau=2,4,\ldots,t-r$. But the value of each of these quotients is
equal to the expression in the right hand side of \eqref{eq:q_2}
with $t$ replaced by $t-\tau$ and $y$ replaced by $y-\frac{\tau}{2}$,
which is strictly smaller than \eqref{eq:q_2}. Hence $\frac{q_{3}}{q_{2}}\le q_{2}$
and $q_{3}=q_{2}\cdot\frac{q_{3}}{q_{2}}\le q_{2}^{~2}$.

A similar argument yields that $\frac{q_{i}}{q_{i-1}}\le q_{2}$ for
all $i\ge3$, hence $q_{i}\le q_{2}^{~i-1}$. Overall we get from
\eqref{eq:precise formula for w_t(o,r)-1} that
\begin{equation}
w_{t}\left(o,r\right)\le\frac{2y-t+1}{y+1}\cdot\frac{d}{d-1}\cdot\binom{t}{y}\left(d-1\right)^{y}\left[1+q_{2}\cdot\frac{d}{d-1}+\left(q_{2}\cdot\frac{d}{d-1}\right)^{2}+\ldots\right].\label{eq:w_t(o,r) with geom series}
\end{equation}
As $d\ge3$ we get from \eqref{eq:q_2 upper bound} that 
\[
q_{2}\cdot\frac{d}{d-1}<\frac{5}{8}\cdot\frac{3}{2}<1.
\]
So the geometric series in \eqref{eq:w_t(o,r) with geom series} is
bounded by a constant, and as $\frac{2y-t+1}{y+1}<1$, there is a
constant $c>0$ with 
\[
w_{t}\left(o,r\right)\le c\binom{t}{y}\left(d-1\right)^{y}.
\]
This completes not only the proof of Lemma \ref{lem:upper-bound, non-constant r},
but also of Proposition \ref{prop:approximate number of walks ending at distance r}.
\end{proof}
We can now prove our main result in its general form -- Theorem \ref{thm:cogrowth for functions}.
We begin with the first item, concerning $\limsup$.
\begin{proof}[Proof of Theorem \ref{thm:cogrowth for functions}, Item \eqref{enu:limsup}]
 Let $0\ne f\colon V\left(\mathbb{T}_{d}\right)\to\mathbb{R}_{\ge0}$
be a non-negative function, not identically zero, with $\alpha\defi\alpha\left(f\right)=\limsup_{t\to\infty}a_{t}\left(f\right)^{1/t}<\infty$.
We ought to show that $\beta\left(f\right)=g_{d}\left(\alpha\right)$. 

We first prove that $\beta(f)\le g_{d}(\alpha)$. If $\alpha<\sqrt{d-1}$,
we may replace $f$ with some $f'$ with $f'-f\ge0$ and $\alpha(f')=\sqrt{d-1}$,
and obviously $\beta(f')\ge\beta(f)$, while $g_{d}(\alpha(f'))=g_{d}(\alpha)=2\sqrt{d-1}$.
So for the sake of proving the upper bound, we may assume without
loss of generality that $\alpha\ge\sqrt{d-1}$. By symmetry, the number
of walks of length $t$ from $o$ to a given vertex $v$, depends
only on $\mathrm{dist}\left(o,v\right)$, so the average value of
$f$ at the end of the walks counted in $w_{t}\left(o,r\right)$ is
$\frac{a_{t}\left(f\right)}{d\left(d-1\right)^{r-1}}$ (or just $a_{0}\left(f\right)$
if $r=0$). Using the function $h\left(t\right)$ from Proposition
\ref{prop:approximate number of walks ending at distance r} and noting
that for every $\varepsilon>0$ there is a constant $c_{2}>0$ with
$a_{r}\left(f\right)\le c_{2}\cdot\left(\alpha+\varepsilon\right)^{r}$
for all $r$, we get
\begin{eqnarray*}
b_{t}\left(f\right) & = & w_{t}\left(o,o\right)a_{0}\left(f\right)+\sum_{r=1}^{t}w_{t}\left(o,r\right)\frac{a_{r}\left(f\right)}{d\left(d-1\right)^{r-1}}\\
 & \stackrel{\mathrm{Prop.}~\ref{prop:approximate number of walks ending at distance r}}{\le} & c_{1}h\left(t\right)\sum_{y=\left\lceil \frac{t}{2}\right\rceil }^{t}\binom{t}{y}\left(d-1\right)^{y}\frac{c_{2}\left(\alpha+\varepsilon\right)^{r}}{\left(d-1\right)^{r}}\\
 & \le & c_{1}c_{2}h\left(t\right)\sum_{y=0}^{t}\binom{t}{y}\left(d-1\right)^{y}\left(\frac{\alpha+\varepsilon}{d-1}\right)^{2y-t}\\
 & = & c_{1}c_{2}h\left(t\right)\sum_{y=0}^{t}\binom{t}{y}\left(\frac{d-1}{\alpha+\varepsilon}\right)^{t-y}\cdot\left(\alpha+\varepsilon\right)^{y}\\
 & = & c_{1}c_{2}h\left(t\right)\left[\frac{d-1}{\alpha+\varepsilon}+\alpha+\varepsilon\right]^{t}=c_{1}c_{2}h\left(t\right)g_{d}\left(\alpha+\varepsilon\right)^{t},
\end{eqnarray*}
where $c_{1}$ is some positive constant. So $\beta\left(f\right)\le g_{d}\left(\alpha+\varepsilon\right)$
for every $\varepsilon>0$, hence by the continuity of $g_{d}$, $\beta\left(f\right)\le g_{d}\left(\alpha\right)$.

Conversely, we prove that $\beta(f)\ge g_{d}(\alpha)$. First consider
the case $\alpha\le\sqrt{d-1}$. For every vertex $v\in V\left(\mathbb{T}_{d}\right)$,
$\beta\left(\mathbb{1}_{\left\{ v\right\} }\right)$ is equal to the
spectral radius of the tree, which is $2\sqrt{d-1}$. As $f\ne0$,
there exists some vertex $v$ with $f\left(v\right)>0$, and then
$\beta\left(f\right)\ge\beta\left(\mathbb{1}_{\left\{ v\right\} }\right)=2\sqrt{d-1}=g_{d}(\alpha)$.

Now assume that $\alpha>\sqrt{d-1}$. There is some subsequence $\left\{ t_{k}\right\} $
with $\lim_{k\to\infty}a_{t_{k}}\left(f\right)^{1/t_{k}}=\alpha$,
so for every $\sqrt{d-1}<\overline{\alpha}<\alpha\left(f\right)$
we have $a_{t_{k}}\left(f\right)>c_{3}\cdot\overline{\alpha}^{t_{k}}$
for some $c_{3}>0$ (and large enough $k$). For every $k$, we find
a positive integer $u_{k}$ so that a 'significant' part of the mass
of $f$ over all length-$u$ walks, is obtained at the length-$u$
walks ending at distance $t_{k}$. Formally, let $u_{k}\in\mathbb{Z}_{\ge0}$
be so that among all values of $\binom{u_{k}}{y}\left(\frac{d-1}{\overline{a}}\right)^{u_{k}-y}\overline{\alpha}^{y}$
where $y=0,1,\ldots,u_{k}$, the value is largest at $y_{k}=\frac{t_{k}+u_{k}}{2}$.\footnote{We remark that for a given $u_{k}$, the largest value of $\binom{u_{k}}{y}\left(d-1\right)^{u_{k}-y}\left(\overline{\alpha}^{2}\right)^{y}$
is obtained for $y_{k}\approx\frac{\overline{\alpha}^{2}}{\overline{\alpha}^{2}+d-1}u_{k}>\frac{1}{2}u_{k}$
(recall that $\overline{\alpha}>\sqrt{d-1}$), so $t_{k}=2y_{k}-u_{k}\approx\varepsilon u_{k}$
with $\varepsilon=\frac{\overline{\alpha}^{2}-\left(d-1\right)}{\overline{\alpha}^{2}+\left(d-1\right)}$.
A simple computation yields that $u_{k}$ can be taken to be any integer
of the same parity as $t_{k}$ in $\left[\frac{1}{\varepsilon}\left(t_{k}-\frac{2\overline{\alpha}^{2}}{\overline{\alpha}^{2}+\left(d-1\right)}\right),\frac{1}{\varepsilon}\left(t_{k}+\frac{2\left(d-1\right)}{\overline{\alpha}^{2}+\left(d-1\right)}\right)\right]$.} 

By the assumption on $u_{k}$ and $y_{k}$ we have\footnote{Up to a constant, one can replace the $\frac{1}{u_{k}+1}$ term in
\eqref{eq:lower bound for maximal val-1} with $\frac{1}{\sqrt{u_{k}}}$,
but we do not need this kind of precision here.}:
\begin{equation}
\binom{u_{k}}{y_{k}}\left(\frac{d-1}{\overline{\alpha}}\right)^{u_{k}-y_{k}}\overline{\alpha}^{y_{k}}\ge\frac{1}{u_{k}+1}\sum_{y=0}^{u_{k}}\binom{u_{k}}{y}\left(\frac{d-1}{\overline{\alpha}}\right)^{u_{k}-y}\overline{\alpha}^{y}=\frac{1}{u_{k}+1}\left(\frac{d-1}{\overline{\alpha}}+\overline{\alpha}\right)^{u_{k}}.\label{eq:lower bound for maximal val-1}
\end{equation}
Therefore, with $c_{4}$ and $c_{5}$ being appropriate constants,
for large enough $k$ we have
\begin{eqnarray}
\beta_{u_{k}}\left(f\right) & \ge & w_{u_{k}}\left(o,t_{k}\right)\frac{a_{t_{k}}\left(f\right)}{d\left(d-1\right)^{t_{k}-1}}\nonumber \\
 & \stackrel{\mathrm{Lemma}~\ref{lem:lower bound}}{\ge} & c_{4}\frac{1}{u_{k}}\binom{u_{k}}{y_{k}}\left(d-1\right)^{y_{k}}\left(\frac{\overline{\alpha}}{d-1}\right)^{t_{k}}=c_{4}\frac{1}{u_{k}}\binom{u_{k}}{y_{k}}\left(d-1\right)^{y_{k}}\left(\frac{\overline{\alpha}}{d-1}\right)^{2y_{k}-u_{k}}\nonumber \\
 & = & c_{4}\frac{1}{u_{k}}\binom{u_{k}}{y_{k}}\left(\frac{d-1}{\overline{\alpha}}\right)^{u_{k}-y_{k}}\overline{\alpha}^{y_{k}}\nonumber \\
 & \stackrel{\eqref{eq:lower bound for maximal val-1}}{\ge} & c_{4}\frac{1}{u_{k}}\cdot\frac{1}{u_{k}+1}\left(\frac{d-1}{\overline{\alpha}}+\overline{\alpha}\right)^{u_{k}}\ge c_{5}\cdot\frac{1}{u_{k}^{~2}}g_{d}\left(\overline{\alpha}\right)^{u_{k}}.\label{eq:lower bound for b_u_k}
\end{eqnarray}
It follows that $\beta\left(f\right)\ge g_{d}\left(\overline{\alpha}\right)$.
Therefore $\beta\left(f\right)\ge g_{d}\left(\alpha\right)$, and
we conclude that $\beta\left(f\right)=g_{d}\left(\alpha\right)$.
\end{proof}

\begin{proof}[Proof of Theorem \ref{thm:cogrowth for functions}, Item \eqref{enu:ordinary limit}]
 Recall that Item \eqref{enu:ordinary limit} in Theorem \ref{thm:cogrowth for functions}
consists of three parts. Note that part \ref{enu:limit-odd-and-even}
is immediate from parts \ref{enu:limit-even} and \ref{enu:limit-odd}.
We prove here part \ref{enu:limit-even}, the proof of part \ref{enu:limit-odd}
being almost identical.

So assume that $\tilde{\alpha}^{\mathrm{even}}\left(f\right)\defi\lim_{t\to\infty}a_{2t}\left(f\right)^{1/(2t)}$
exists and that $f|_{\mathbb{T}_{d}^{\mathrm{even}}}\ne0$. By the
proof of Item \eqref{enu:limsup} of the theorem, $\limsup_{t\to\infty}b_{2t}\left(f\right)^{1/\left(2t\right)}\le g_{d}\left(\tilde{\alpha}^{\mathrm{even}}\left(f\right)\right)$
(note that the end vertex of an even-length walk is at an even distance
from $o$). Let $v\in\mathbb{T}_{d}^{\mathrm{even}}$ be a vertex
with $f\left(v\right)>0$. Then $\liminf_{t\to\infty}b_{2t}\left(f\right)^{1/\left(2t\right)}\ge\lim_{t\to\infty}b_{2t}\left(\mathbb{1}_{\left\{ v\right\} }\right)^{1/\left(2t\right)}=\rho\left(\mathbb{T}_{d}\right)=2\sqrt{d-1}$.
So if $\tilde{\alpha}^{\mathrm{even}}\left(f\right)\le\sqrt{d-1}$,
we established that $\tilde{\beta}^{\mathrm{even}}\left(f\right)=\lim_{t\to\infty}\beta_{2t}\left(f\right)^{1/\left(2t\right)}$
exists and equals $2\sqrt{d-1}$. 

Now assume that $\tilde{\alpha}^{\mathrm{even}}\left(f\right)>\sqrt{d-1}$.
We may repeat the argument in the proof of Item \eqref{enu:limsup}
with $u_{k}=2k$ and arbitrary $\sqrt{d-1}<\overline{\alpha}<\tilde{\alpha}^{\mathrm{even}}\left(f\right)$.
As above, $y_{k}\approx\frac{\overline{\alpha}^{2}}{\overline{\alpha}^{2}+d-1}\cdot2k$
is the value maximizing $\binom{2k}{y}\left(\frac{d-1}{\overline{\alpha}}\right)^{2k-y}\overline{\alpha}^{y}$,
and then the same argument from \eqref{eq:lower bound for b_u_k}
shows that $\beta_{2k}\left(f\right)$ is at least a constant times
$\frac{1}{k^{2}}g_{d}\left(\overline{\alpha}\right)^{2k}$. We conclude
that $\liminf_{t\to\infty}\beta_{2t}\left(f\right)^{1/\left(2t\right)}\ge g_{d}\left(\tilde{\alpha}^{\mathrm{even}}\left(f\right)\right)$.
Therefore, $\tilde{\beta}^{\mathrm{even}}\left(f\right)$ exists and
equals $g_{d}\left(\tilde{\alpha}^{\mathrm{even}}\left(f\right)\right)$.
\end{proof}

\section{The cogrowth formula for arbitrary subsets of the biregular tree\label{sec:biregular}}

In this section we prove Theorem \ref{thm:cogrowth biregular} which
extends Theorem \ref{thm:cogrowth for subsets} to subsets of the
biregular tree. Let $2\le c\le d$ be integers with $3\le d$ (when
$3\le c=d$ this recovers the $d$-regular case). Recall that $\mathbb{T}_{c,d}$
is the $\left(c,d\right)$-biregular with a root vertex $o$. To shorten
the formulas, we use the notation 
\[
\overline{c}\defi c-1~~~~~~~\mathrm{and}~~~~~~~\overline{d}\defi d-1.
\]
Define $g_{c,d}\colon\mathbb{R}_{\ge0}\to\mathbb{R}_{\ge\sqrt{\overline{c}}+\sqrt{\overline{d}}}$
by 
\begin{equation}
g_{c,d}\left(\alpha\right)=\begin{cases}
\sqrt{\overline{c}}+\sqrt{\overline{d}} & \mathrm{if}~\alpha\le\left(\overline{c}\cdot\overline{d}\right)^{1/4},\\
\sqrt{\alpha{}^{2}+\overline{c}+\overline{d}+\frac{\overline{c}\cdot\overline{d}}{\alpha{}^{2}}} & \mathrm{if}~\alpha\ge\left(\overline{c}\cdot\overline{d}\right)^{1/4}.
\end{cases}\label{eq:def of g_c,d}
\end{equation}
Our goal is to prove that for any non-empty subset $\emptyset\ne S\subseteq V\left(\mathbb{T}_{c,d}\right)$
of the vertices of $\mathbb{T}_{c,d}$, we have $\beta\left(S\right)=g_{c,d}\left(\alpha\left(S\right)\right)$. 

We believe that the more general result concerning non-negative functions,
as in Theorem \ref{thm:cogrowth for functions}, should hold, but
we were not able to make our arguments work in this generality\footnote{The part we were not able to generalize to the biregular case is the
analog of Lemma \ref{lem:upper-bound, non-constant r}.}. 
\begin{conjecture}
Fix $2\le c\le d$, $3\le d$ and let $0\ne f\colon V\left(\mathbb{T}_{c,d}\right)\to\mathbb{R}_{\ge0}$
be a non-negative function on the vertices of the $\left(c,d\right)$-biregular
tree, not identically zero. Then $\beta\left(f\right)=g_{c,d}\left(\alpha\left(f\right)\right)$.
\end{conjecture}

The fact we restrict attention to subsets of the tree (namely, to
indicator function rather than arbitrary non-negative functions) carries
the advantage of significantly simplifying our proof, as compared
with the proof in Section \ref{sec:regular}. We will point out the
places in the proof where we need this restriction.

We grade $\mathbb{T}_{c,d}$ similarly to the grading of $\mathbb{T}_{d}$:
the root $o$ is at level $o$, and every vertex $v$ has a single
parent sitting one level down, and $\deg\left(v\right)-1$ children
one level up. Our analysis of walks is of two steps, or a ``double-step'',
at a time. There is a total of $cd$ options for every double-step.
Of these, the level changes by
\[
\begin{cases}
-2 & \mathrm{in}~1~~\mathrm{option}\\
0 & \mathrm{in}~\overline{c}+\overline{d}~\mathrm{~options}\\
2 & \mathrm{in}~\overline{c}\cdot\overline{d}~~\mathrm{options.}
\end{cases}
\]
Suppose that in a walk of length $2t$, out of $t$ double-steps,
there are $x$ downward double-steps, $y$ ``horizontal'' double-steps
and $z$ upward double-steps. The level at the end is then $2z-2x$.
As in the regular case, we say that a walk of length $2t$ is\emph{
}\textbf{positive} if it never returns to level zero at a positive
time. 
\begin{lem}
\label{lem:strictly positive biregular}For every $x,y,z\in\mathbb{Z}_{\ge0}$
with $z-x>0$, out of the $\binom{x+y+z}{x~~y~~z}\left(\overline{c}+\overline{d}\right)^{y}\left(\overline{c}\cdot\overline{d}\right)^{z}$
walks in $\mathbb{T}_{c,d}$ with precisely $x$ downward double-steps,
$y$ horizontal double-steps and $z$ upward double-steps, the proportion
of positive walks is 
\[
\frac{z-x}{x+y+z}.
\]
\end{lem}

\begin{proof}
Denote $t=x+y+z$. The argument is the same as in Lemma \ref{lem:positive and non-negative walks}:
An even-length positive walk begins with an upward double-step, and
then every double step does not end at height below $2$. Among all
walks with first double-step upward, those that are not positive are
in bijection with walks from level $2$ to level $2x-2z$ with $t-1$
double-steps and $y$ horizontal double steps among them. So the proportion
of positive walks is
\[
\frac{\binom{t-1}{x~y~z-1}-\binom{t-1}{z~y~x-1}}{\binom{t}{x~y~z}}=\frac{z-x}{t}.
\]
\end{proof}
\begin{proof}[Proof of Theorem \ref{thm:cogrowth biregular}]
 Let $\emptyset\ne S\subseteq V\left(\mathbb{T}_{c,d}\right)$ be
an arbitrary non-empty subset of the vertices of $\mathbb{T}_{c,d}$.
Denote $\alpha\defi\alpha\left(S\right)$. The maximal possible value
of $\alpha$ is obtained when $S=V(\mathbb{T}_{c,d})$, in which case
$\alpha=\sqrt{\overline{c}\cdot\overline{d}}$. So, in general, $\alpha\in\left[0,\sqrt{\overline{c}\cdot\overline{d}}\right]$.
We need to show that $\beta\left(S\right)=g_{c,d}\left(\alpha\right)$. 

First we show that $\beta\left(S\right)\le g_{c,d}\left(\alpha\right)$.
If $\alpha=\sqrt{\overline{c}\cdot\overline{d}}$ has the maximal
possible value, then $g_{c,d}\left(\alpha\right)=\sqrt{cd}=\beta\left(V\left(\mathbb{T}_{c,d}\right)\right)$
is the maximal possible value of $\beta$ on any set $S$, and, obviously,
$\beta\left(S\right)\le g_{c,d}\left(\alpha\right)$.\footnote{This is one minor spot where the restriction to indicator functions
is used.} Now assume that $\alpha<\sqrt{\overline{c}\cdot\overline{d}}$. Let
$\overline{\alpha}$ be an arbitrary number with $\left(\overline{c}\cdot\overline{d}\right)^{1/4}<\overline{\alpha}<\left(\overline{c}\cdot\overline{d}\right)^{1/2}$
and $\alpha<\overline{\alpha}$. There exists a constant $k_{\overline{\alpha}}>0$
such that for all $t$ we have $a_{t}\left(S\right)\le k_{\overline{\alpha}}\overline{\alpha}^{t}$.
By symmetry, instead of considering $S$ or its indicator function,
we may consider the spherical function $f\colon V\left(\mathbb{T}_{c,d}\right)\to\mathbb{R}_{\ge0}$
defined by
\[
f\left(v\right)=\frac{a_{t}\left(S\right)}{\#\left\{ \text{vertices~in~}\mathbb{T}_{c,d}\,\,\mathrm{at~distance}~t~\mathrm{from~}o\right\} }
\]
for every vertex $v$ at (arbitrary) distance $t$ from $o$. Let
us start with walks of even length. For a walk $w$ (always starting
at $o$), we denote the distance of its end vertex from $o$ by $\delta\left(w\right)$,
and the level of the end vertex in the graded tree by $\lambda\left(w\right)$.
Clearly, $\lambda(w)\le\delta(w)$. Denote by $\lessapprox$ an inequality
up to a universal constant. For every $t\in\mathbb{Z}_{\ge1}$,

\begin{eqnarray}
b_{2t}\left(S\right)=b_{2t}\left(f\right) & = & \sum_{\substack{w:\,\mathrm{walk\,of}\\
\mathrm{length}~2t
}
}f\left(\text{end-vertex-of}\,w\right)\nonumber \\
 & \lessapprox & k_{\overline{\alpha}}\sum_{\substack{w:\,\mathrm{walk\,of}\\
\mathrm{length}~2t
}
}\left(\frac{\overline{\alpha}}{\sqrt{\overline{c}\cdot\overline{d}}}\right)^{\delta\left(w\right)}\nonumber \\
 & \stackrel{\overline{\alpha}<\sqrt{\overline{c}\cdot\overline{d}},\lambda\left(w\right)\le\delta\left(w\right)}{\le} & k_{\overline{\alpha}}\sum_{\substack{w:\,\mathrm{walk\,of}\\
\mathrm{length}~2t
}
}\left(\frac{\overline{\alpha}}{\sqrt{\overline{c}\cdot\overline{d}}}\right)^{\lambda\left(w\right)}\label{eq:where we use restriction to indicator function}\\
 & = & k_{\overline{\alpha}}\sum_{\substack{x+y+z=t\\
x,y,z\ge0
}
}\binom{t}{x~y~z}\left[\overline{c}+\overline{d}\right]^{y}\left[\overline{c}\cdot\overline{d}\right]^{z}\left(\frac{\overline{\alpha}}{\sqrt{\overline{c}\cdot\overline{d}}}\right)^{2z-2x}\nonumber \\
 & = & k_{\overline{\alpha}}\sum_{\substack{x+y+z=t\\
x,y,z\ge0
}
}\binom{t}{x~y~z}\left[\frac{\overline{c}\cdot\overline{d}}{\overline{\alpha}^{2}}\right]^{x}\left(\overline{c}+\overline{d}\right)^{y}\cdot\left(\overline{\alpha}^{2}\right)^{z}\label{eq:cogrowth-upper-bound-1}\\
 & = & k_{\overline{\alpha}}\left[\frac{\overline{c}\cdot\overline{d}}{\overline{\alpha}^{2}}+\overline{c}+\overline{d}+\overline{\alpha}^{2}\right]^{t}=k_{\overline{\alpha}}g_{c,d}\left(\overline{\alpha}\right)^{2t}.\nonumber 
\end{eqnarray}
(The inequality \eqref{eq:where we use restriction to indicator function}
is the main point where we use out restriction to indicator functions,
which guarantess that $\alpha\le\sqrt{\overline{c}\cdot\overline{d}}$.)
The same argument works for odd-length talks: by symmetry, we may
assume the first step of a length-$\left(2t+1\right)$ walk is deterministically
upward, and then break the rest of the walk into double-steps as before.
We get the same computation up to a constant, so $b_{2t+1}\left(f\right)\lessapprox k_{\overline{\alpha}}g_{c,d}\left(\overline{\alpha}\right)^{2t+1}$.
This proves that $\beta\left(S\right)\le g_{c,d}\left(\overline{\alpha}\right)$
for all $\overline{\alpha}$ as above, so $\beta\left(S\right)\le g_{c,d}\left(\alpha\right)$.

Next, we prove that $\beta\left(S\right)\ge g_{c,d}\left(\alpha\right)$.
First, for every $v\in S$, $\beta\left(\mathbb{1}_{\left\{ v\right\} }\right)$
is equal to the spectral radius of the tree, which is $\sqrt{\overline{c}}+\sqrt{\overline{d}}$
\cite{GODSIL1988}, so $\beta\left(S\right)\ge\sqrt{\overline{c}}+\sqrt{\overline{d}}$.
This proves the lower bound whenever $\alpha\le\left(\overline{c}\cdot\overline{d}\right)^{1/4}$. 

Assume from now on that $\left(\overline{c}\cdot\overline{d}\right)^{1/4}<\alpha$.
Assume without loss of generality that there is a subsequence of \uline{even}
positive integers $\left\{ 2t_{k}\right\} _{k=1}^{\infty}$ with $\alpha=\lim_{k\to\infty}a_{2t_{k}}\left(S\right)$.
Let $\overline{\alpha}$ be arbitrary in $\left(\left(\overline{c}\cdot\overline{d}\right)^{1/4},\alpha\right)$.
For every large enough $k$, $a_{2t_{k}}\left(f\right)\ge\overline{\alpha}^{2t_{k}}$. 

Let $u_{k}$ be a positive integer so that length-$2u_{k}$ walks
on the graded tree $\mathbb{T}_{c,d}$ are relatively likely to end
at level $2t_{k}$. Formally, $u_{k}$ is a positive integer so that
among all triples of non-negative integers $x,y,z$ with $x+y+z=u_{k}$,
the quantity
\[
\binom{u_{k}}{x~y~z}\left[\frac{\overline{c}\cdot\overline{d}}{\overline{\alpha}^{2}}\right]^{x}\left(\overline{c}+\overline{d}\right)^{y}\cdot\left(\overline{\alpha}^{2}\right)^{z}
\]
is maximized at $x_{k},y_{k},z_{k}$ satisfying $2z_{k}-2x_{k}=2t_{k}$.
Generally, $\binom{u}{x\,y\,z}p^{x}q^{y}r^{z}$ is maximized when
$x:y:z$ is roughly $p:q:r$. In our case, we have that $u_{k}$,
$x_{k}$, $y_{k}$ and $z_{k}$ are roughly $\frac{\overline{c}\cdot\overline{d}+\overline{\alpha}^{2}\left(\overline{c}+\overline{d}\right)+\overline{\alpha}^{4}}{\overline{\alpha}^{4}-\overline{c}\cdot\overline{d}}t_{k}$,
$\frac{\overline{c}\cdot\overline{d}}{\overline{\alpha}^{4}-\overline{c}\cdot\overline{d}}t_{k}$,
$\frac{\overline{\alpha}^{2}(\overline{c}+\overline{d})}{\overline{\alpha}^{4}-\overline{c}\cdot\overline{d}}t_{k}$
and $\frac{\overline{\alpha}^{4}}{\overline{\alpha}^{4}-\overline{c}\cdot\overline{d}}t_{k}$,
respectively. We obtain 
\begin{equation}
\binom{u_{k}}{x_{k}~y_{k}~z_{k}}\left[\frac{\overline{c}\cdot\overline{d}}{\overline{\alpha}^{2}}\right]^{x_{k}}\left(\overline{c}+\overline{d}\right)^{y_{k}}\cdot\left(\overline{\alpha}^{2}\right)^{z_{k}}\gtrapprox\frac{1}{u_{k}^{~2}}\left[\frac{\overline{c}\cdot\overline{d}}{\overline{\alpha}^{2}}+\overline{c}+\overline{d}+\overline{\alpha}^{2}\right]^{u_{k}}.\label{eq:maximal y_k,z_k}
\end{equation}

By Lemma \ref{lem:strictly positive biregular}, out of all length-2$u_{k}$
walks with $x_{k}$, $y_{k}$ and $z_{k}$ double steps that are downward,
horizontal and upward, respectively, the proportion of positive walks
is 
\begin{equation}
\frac{z_{k}-x_{k}}{u_{k}}=\frac{t_{k}}{u_{k}}>\frac{1}{u_{k}},\label{eq:proportion of positive}
\end{equation}
the inequality holding for almost every $k$. Now, for every large
enough $k$,

\begin{eqnarray*}
b_{2u_{k}}\left(f\right) & = & \sum_{\substack{w:\,\mathrm{walk\,of}\\
\mathrm{length}\,2u_{k}
}
}f\left(\text{end-vertex-of}\,w\right)~\\
 & \ge & \sum_{\substack{w:\,\mathrm{walk~of~length}~2u_{k}~\mathrm{with~exactly}\\
x_{k},y_{k},z_{k}~\mathrm{downard/horizontal/upward~double-steps}
}
}f\left(\text{end-vertex-of}\,w\right)\\
 & \ge & \sum_{\substack{w:\,\mathrm{positive~walk~of~length}~2u_{k}~\mathrm{with~exactly}\\
x_{k},y_{k},z_{k}~\mathrm{downard/horizontal/upward~double-steps}
}
}f\left(\text{vertex~at~distance~}2z_{k}-2x_{k}\right)\\
 & \stackrel{\eqref{eq:proportion of positive}}{\gtrapprox} & \frac{1}{u_{k}}\binom{u_{k}}{x_{k}~y_{k}~z_{k}}\left(\overline{c}+\overline{d}\right)^{y_{k}}\left[\overline{c}\cdot\overline{d}\right]^{z_{k}}\left(\frac{\overline{\alpha}^{2}}{\overline{c}\cdot\overline{d}}\right)^{z_{k}-x_{k}}\\
 & = & \frac{1}{u_{k}}\binom{u_{k}}{x_{k}~y_{k}~z_{k}}\left[\frac{\overline{c}\cdot\overline{d}}{\overline{\alpha}^{2}}\right]^{x_{k}}\left(\overline{c}+\overline{d}\right)^{y_{k}}\left(\overline{\alpha}^{2}\right)^{z_{k}}\\
 & \stackrel{\eqref{eq:maximal y_k,z_k}}{\ge} & \frac{1}{u_{k}^{~3}}\left[\frac{\overline{c}\cdot\overline{d}}{\overline{\alpha}^{2}}+\overline{c}+\overline{d}+\overline{\alpha}^{2}\right]^{u_{k}}=\frac{1}{u_{k}^{~3}}g_{c,d}\left(\overline{\alpha}\right)^{2u_{k}}.
\end{eqnarray*}
Hence $\beta\left(S\right)\ge g_{c,d}\left(\overline{\alpha}\right)$
for arbitrary $\overline{\alpha}\in\left(\left(\overline{c}\cdot\overline{d}\right)^{1/4},\alpha\right)$,
and we conclude that $\beta\left(S\right)\ge g_{c,d}\left(\alpha\right)$.
\end{proof}
\bibliographystyle{amsalpha}
\bibliography{united,self}

\end{document}